\documentclass[11pt]{article}

\usepackage{hyperref}
\usepackage{verbatim}
\usepackage{xcolor}
\usepackage{graphicx}
\usepackage{float}
\usepackage{lmodern}
\usepackage{amsthm}
\usepackage{cite}
\usepackage{amsmath}
\usepackage{mathtools}
\usepackage{amssymb}
\usepackage{bbm}
\usepackage{subfig}
\usepackage{wrapfig}
\usepackage{sepfootnotes}

\DeclarePairedDelimiter\bra{\langle}{\rvert}
\DeclarePairedDelimiter\ket{\lvert}{\rangle}
\DeclarePairedDelimiterX\braket[2]{\langle}{\rangle}{#1 \delimsize\vert #2}

\hypersetup{
pdfpagemode={UseOutlines},
bookmarksopen,
pdfstartview={FitH},
colorlinks,
linkcolor={blue},
citecolor={red},
urlcolor={blue}
}

\graphicspath{{./figures/}}

\numberwithin{equation}{section}

\theoremstyle{plain} 

\newtheorem{theorem}{Theorem}[section] 
\newtheorem{lemma}[theorem]{Lemma} 
\newtheorem{remark}[theorem]{Remark}

\newcommand{\GG}{\mathbb{G}}
\newcommand{\Pro}{\mathbb{P}}

\newcommand{\Exp}{\mathbb{E}}
\newcommand{\Var}{\mathrm{Var}}
\newcommand{\1}{\mathbbm{1}}
\newcommand{\ones}{\vec{\1}}
\newcommand{\ER}{Erd\H{o}s-R\'enyi\ }
\newcommand{\lep}{\left(}
\newcommand{\rip}{\right)}
\newcommand{\Real}{\mathbb{R}}

\newcommand{\Ecan}{\Exp_{\text{can}}}
\newcommand{\Emic}{\Exp_{\text{mic}}}
\newcommand{\Pmic}{\Pro_{\text{mic}}}
\newcommand{\Pcan}{\Pro_{\text{can}}}

\newcommand{\ee}{\mathrm{e}}

\newcommand\blfootnote[1]{%
  \begingroup
  \renewcommand\thefootnote{}\footnote{#1}%
  \addtocounter{footnote}{-1}%
  \endgroup
}

\DeclarePairedDelimiter\abs{\lvert}{\rvert}
\DeclarePairedDelimiter\norm{\lVert}{\rVert}

\makeatletter
\let\oldabs\abs
\def\abs{\@ifstar{\oldabs}{\oldabs*}}

\let\oldnorm\norm
\def\norm{\@ifstar{\oldnorm}{\oldnorm*}}
\makeatother

\begin{document}

\title{\vspace{-55pt}A spectral signature of breaking of ensemble equivalence\\ 
for constrained random graphs}

\author{
\renewcommand{\thefootnote}{\arabic{footnote}}
P. Dionigi
\footnotemark[1]
\\
\renewcommand{\thefootnote}{\arabic{footnote}}
D. Garlaschelli
\footnotemark[2] 
\\
\renewcommand{\thefootnote}{\arabic{footnote}}
F. den Hollander
\footnotemark[1] 
\\
\renewcommand{\thefootnote}{\arabic{footnote}}
M. Mandjes
\footnotemark[4]
\blfootnote{E-mails: \href{mailto:p.dionigi@math.leidenuniv.nl}{p.dionigi@math.leidenuniv.nl}; \href{mailto:garlaschelli@lorentz.leidenuniv.nl}{garlaschelli@lorentz.leidenuniv.nl}; \href{mailto:denholla@math.leidenuniv.nl}{denholla@math.leidenuniv.nl}; \href{mailto:M.R.H.Mandjes@uva.nl}{M.R.H.Mandjes@uva.nl}}
}

\footnotetext[1]{
Mathematical Institute, Leiden University,\\ 
P.O.\ Box 9512, 2300 RA Leiden, The Netherlands 
}

\footnotetext[2]{
Lorentz Institute for Theoretical Physics, Leiden University,\\ 
P.O.\  Box 9504, 2300 RA Leiden, The Netherlands
}

\footnotetext[3]{
IMT School for Advanced Studies,\\
Piazza S.\ Francesco 19, 55100 Lucca, Italy
}

\footnotetext[4]{
Korteweg-de Vries Institute, University of Amsterdam,\\ 
P.O.\ Box 94248, 1090 GE Amsterdam, The Netherlands
}

\date{8 July 2021}

\maketitle
\vspace{-20pt}
\begin{abstract}
\noindent
For random systems subject to a constraint, the \emph{microcanonical ensemble} requires the constraint to be met by every realisation (`hard constraint'), while the \emph{canonical ensemble} requires the constraint to be met only on average (`soft constraint'). It is known that for random graphs subject to topological constraints \emph{breaking of ensemble equivalence} may occur when the size of the graph tends to infinity, signalled by a non-vanishing \emph{specific relative entropy} of the two ensembles. We investigate to what extent breaking of ensemble equivalence is manifested through the \emph{largest eigenvalue} of the adjacency matrix of the graph. We consider two examples of constraints in the dense regime: (1) fix the degrees of the vertices (= the degree sequence); (2) fix the sum of the degrees of the vertices (= twice the number of edges). Example (1) imposes an extensive number of local constraints and is known to lead to breaking of ensemble equivalence. Example (2) imposes a single global constraint and is known to lead to ensemble equivalence. Our working hypothesis is that breaking of ensemble equivalence corresponds to a \emph{non-vanishing difference of the expected values of the largest eigenvalue} under the two ensembles. We verify that, in the limit as the size of the graph tends to infinity, the difference between the expected values of the largest eigenvalue in the two ensembles does not vanish for (1) and vanishes for (2). A key tool in our analysis is a transfer method that uses relative entropy to determine whether probabilistic estimates can be carried over from the canonical ensemble to the microcanonical ensemble, and illustrates how breaking of ensemble equivalence may prevent this from being possible.  \vspace{0.2cm}\\
\noindent
\emph{Key words.}\\
\small Constrained random graphs; canonical and micro-canonical ensembles; ensemble equivalence; relative entropy; adjacency matrix; largest eigenvalue; Erd\H{o}s-R\'enyi random graph; regular random graph.\\
\emph{MSC2020:}
60B20, 
60C05, 
60K35. 
\\
\emph{Acknowledgment.}
The work of PD, FdH and MM was supported by the Netherlands Organisation for Scientific Research (NWO) through Gravitation-grant NETWORKS-024.002.003, the work of DG by the Dutch Econophysics Foundation (Stichting Econophysics, Leiden, The Netherlands).
\normalsize
\end{abstract}

\newpage

\section{Introduction}
\label{sec:intro}

\paragraph{Background.}

Spectral properties of random graphs have been studied intensively in past years. A non-exhaustive list of key contributions is \cite{erdosSpectralStatisticsErdosRenyi2012,erdosSpectralStatisticsErdos2013,alt_extremal_2020,benaych-georges_largest_2017,chakrabarty2019,dhara_large_2020,ding2010spectral,furedi_eigenvalues_1981,krivelevich_largest_2003, zhu2020}. Both the adjacency matrix and the Laplacian matrix have been popular. Scaling properties have been derived for the spectral distribution and the largest eigenvalue, with focus on central limit and large deviation behaviour. Most papers deal with random graphs whose edges are drawn \emph{independently}. Different types of behaviour show up in the \emph{dense regime} (where the number of edges is of the order of the square of the number of vertices), in the \emph{sparse regime} (where the number of edges is of the order of the number of vertices), and in between. 

In this paper we focus on the \emph{largest eigenvalue} of the non-normalized and non-centred adjacency matrix for a class of \emph{constrained} random graphs. The largest eigenvalue is a highly non-linear functional of the entries of the adjacency matrix and therefore carries global information about the structure of the graph. Constraints are natural in the framework of statistical mechanics and \emph{Gibbs ensembles}. Typically, they introduce a dependence between the edges that makes the spectral analysis challenging.


\paragraph{Breaking of ensemble equivalence (BEE).}

One of the interesting phenomena exhibited by certain classes of constrained random graphs is \emph{Breaking of Ensemble Equivalence} (BEE). To understand what this is, we recall that in statistical physics different microscopic descriptions are available for a system that is subjected to a constraint, referred to as \emph{Gibbs ensembles}. In the \emph{microcanonical ensemble} the constraint is \emph{hard}, i.e., each microscopic realisation of the system matches the constraint \emph{exactly}. In the \emph{canonical ensemble} the constraint is \emph{soft}, i.e., is met only \emph{on average}. For finite systems the two ensembles are clearly different, since they represent different physical situations (energetic isolation, respectively, thermal equilibrium with a reservoir at an appropriate temperture). However, the general belief is that this discrepancy vanishes in the thermodynamic limit. This expectation, referred to as \emph{Equivalence of Ensembles} (EE), permeates the theory of Gibbs ensembles. It turns out that for many physical systems EE holds, but not for all. We refer to \cite{touchette_equivalence_2015} for more background.  

For interacting particle systems, EE has been studied at three different levels: thermodynamic, macrostate and measure. It was shown in \cite{touchette_equivalence_2015} that these levels are equivalent. The present paper uses the measure level, which is based on the vanishing of the specific relative entropy. In \cite{garlaschelli_covariance_2018,squartini_breaking_2015,starreveld_2018,squartini_garlaschelli_2020}, the phenomenon of BEE was studied for random graphs subject to different types of constraints. It was found that, interestingly, BEE is the rule rather than the exception for constraints that are either \emph{extensive} in the number of vertices or \emph{frustrated}. An overview can be found in \cite{roccaverde_breaking_2018}.


\paragraph{Spectral signature of BEE.}

Let $A$ be the adjacency matrix of a random graph on $n$ vertices, i.e., $A=\{a_{ij}\}_{i,j \in [n]}$ with $a_{ij} = 1_{\{i \sim j\}}$. Let $\lambda_1(n)$ denote its largest eigenvalue. For $i \in [n]$, let $k_i$ be the degree of vertex $i$. Write $\Exp_{\mathrm{can}}$ and $\Exp_{\mathrm{mic}}$ to denote expectation with respect to the canonical, respectively, microcanonical ensemble. Put
\begin{equation}
\Delta_\infty = \lim_{n\to\infty} \Big(\Exp_{\mathrm{can}}[\lambda_1(n)]-\Exp_{\mathrm{mic}}[\lambda_1(n)]\Big).
\end{equation}
Our \emph{working hypothesis} is that
\begin{equation}
\begin{array}{llll}
&\Delta_\infty \neq 0 &\Longrightarrow &\mathrm{BEE},\\ 
&\mathrm{BEE} &\Longrightarrow &\Delta_\infty \neq 0 \quad \text{ apart from exceptional constraints}.
\end{array}
\end{equation}
The goal of the present paper and future work is to verify when this working hypothesis is valid and to identify what are the exceptional constraints (see Remark~\ref{rmk:BEEregimes} below).
 
We will verify the working hypothesis for two specific examples of constraints: (1) fix the degrees of the vertices (= the degree sequence); (2) fix the sum of the degrees of the vertices (= twice the number of edges). Example (1) corresponds to the so-called \emph{configuration model}. We consider the particular case where all the degrees are fixed at a common value $d(n)$, in which case the microcanonical ensemble becomes the $d(n)$-regular random graph, for which $\lambda_1(n)=d(n)$ with probability 1. For this case, BEE is known to occur for all choices of $d(n) \neq \{0,n-1\}$, and we will see that $\Delta_\infty \neq 0$ except in the ultra-dense regime where $\lim_{n\to\infty} n^{-1}d(n) =1$. The failure of our working hypothesis in this regime is a consequence of the saturation of the adjacency matrix. Indeed, the largest eigenvalue becomes ineffective in detecting BEE when the two ensembles asymptotically concentrate around the complete graph, for which the largest eigenvalue achieves the maximal value $n-1$. In contrast, relative entropy manages to detect BEE because the two ensembles still look different in the ultra-dense regime, where the number of achievable graphs scales as the exponential of $n^2$. For Example (2) we will see that no BEE occurs and that $\Delta_\infty = 0$. For both examples the canonical ensemble coincides with the \ER random graph with an appropriate retention probability~\cite{squartini_breaking_2015}.

For \ER random graphs, $\lambda_1(n)$ was studied for various different regimes in \cite{furedi_eigenvalues_1981,erdosSpectralStatisticsErdos2013,krivelevich_largest_2003}. Throughout the sequel we consider the regime  
\begin{equation}
\label{regime}
\exists\, \beta \in (6,\infty)\colon \quad  
n^{-1} (\log n)^{\beta} \leq p(n)<1-n^{-1} (\log n)^{\beta}.
\end{equation} 

\begin{theorem} {{\rm \cite[Theorem 6.2]{erdosSpectralStatisticsErdos2013}}}
\label{thm:furedi}
Let $G(n,p(n))$ be the \ER random graph on $n$ vertices with retention probability $p(n)$ satisfying \eqref{regime}. Let $\lambda_1(n)$ be the largest eigenvalue of the adjacency matrix of $G(n,p(n))$. Then
\begin{equation}
\label{CLT}
\Exp_{G(n,p(n))}[\lambda_1(n)] = (n-1)p(n) + (1-p(n)) + O\lep\frac{(1-p(n))^{3/2}}{q(n)\sqrt{(n-1)p(n)}}\rip,
\end{equation} 
where $q(n)=\sqrt{(n-1)p(n)}$ when $p(n) \leq c<1$, and $q(n)=\sqrt{(n-1)(1-p(n))}$ when $p(n)=1-o(1)$.
\end{theorem}

\noindent
To state \eqref{CLT}, we removed diagonal entries so as to get simple graphs, as explained in Chapter \ref{sec:erdos}. Theorem \ref{thm:furedi} shows that the largest eigenvalue of the \ER random graph is a perturbative correction around the mean degree $d(n)=(n-1)p(n)$. In the dense regime $p(n) \equiv p \in (0,1)$ we get the classical result from \cite{furedi_eigenvalues_1981}. In the ultra-dense regime, where the complementary graph is sparse, we can still use \cite[Definition 2.1]{erdosSpectralStatisticsErdos2013}. The lower bound on $p(n)$ in \eqref{regime} implies that we do not capture the sparse regime below the connectivity threshold: a crossover in the scaling behaviour of $\lambda_1(n)$ occurs when $d(n) \asymp \log n$, as proved in \cite{alt_extremal_2020}. 

Theorem \ref{thm:furedi} leads us to our main result.

\begin{theorem}
\label{BEEspec}
Let $p(n)$ satisfy \eqref{regime}.\\
{\rm (1)} Let the constraint be $k_i=d(n)$, $i \in [n]$, with $nd(n)$ even and $\lim_{n\to\infty} [n^{-1}d(n)]/p(n)=1$. Then
\begin{equation}
\label{n1}
\Delta_\infty = \left\{\begin{array}{ll}
1-p,
&\textnormal{if } p(n) \equiv p \in(0,1),\\
1, 
&\textnormal{if } p(n)=o(1),\\
0,
&\textnormal{if } p(n)=1-o(1).
\end{array}
\right.
\end{equation}
{\rm (2)} Let the constraint be $\tfrac12\sum_{i \in [n]} k_i=L(n)$ with $\lim_{n\to\infty} [2n^{-2}L(n)]/p(n)=1$. Then
\begin{equation}
\label{n2}
\Delta_\infty = 0.
\end{equation}
\end{theorem}

\noindent
The restriction that $nd(n)$ is even is needed to make the constraint \emph{graphical}, i.e., there exist simple graphs that meet the constraint. Note the remarkable fact that both $\Exp_{\rm mic}[\lambda_1(n)]$ and $\Exp_{\rm can}[\lambda_1(n)]$ tend to infinity as $n\to\infty$ while their difference remains bounded. 

As shown in \cite{garlaschelli_covariance_2018,squartini_breaking_2015}, BEE occurs in example (1) and EE in example (2), and hence Theorem~\ref{BEEspec} supports our working hypothesis that BEE corresponds to a non-vanishing difference of the expected largest eigenvalues under the two ensembles. 

\begin{remark}
\label{rmk:regimes}
In \cite{krivelevich_largest_2003} a general technique is used that also covers the regime $0<p(n)<n^{-1}(\log n)^\beta$. However, as stated by the authors in their conclusions, their method does not allow for a derivation of the asymptotics of $\Exp[\lambda_1(n)]$. Nevertheless, it is worth mentioning that when $p(n)=\frac{c}{n}$, $c \in (0,\infty)$, the asymptotic behaviour of $\lambda_1(n)$ in the \ER model $G(n,p(n))$ is 
\begin{equation}
\label{eq:sublog}
\lim_{n\to\infty} \lep\lambda_1(n)-\sqrt{\frac{\log n}{\log\log n}}\,\,\rip = 0
\end{equation}
with high probability. Interestingly, in view of the results in Section~\ref{sec:BEE}, this suggests that \eqref{n1} may have limit $\infty$ in this regime.
\end{remark}

\begin{remark}
\label{rmk:BEEregimes}
In \cite{garlaschelli_covariance_2018} it is shown that BEE occurs for three regimes of constant degree $d(n)$: (I) $d(n) = o(\sqrt{n}\,)$ (sparse regime); (II) $\delta n \leq d(n) \leq 1-\delta$ for some $\delta \in (0,\tfrac12]$ (dense regime); (III) $d(n) = n-o(\sqrt{n}\,)$ (ultra dense regime). The scaling of the specific relative entropy is $n$ for regimes (I) and (II), and $n \log n$ for regime (II). Theorem~\ref{BEEspec}(1) shows that our working hypothesis holds in regime (I) (subject to $d(n) \geq (\log n)^{\beta}$) and (II), but fails in regime (III). The reason is that, while the specific relative entropy is invariant under the map where edges are replaced by non-edges and vice versa, the same is not true for the largest eigenvalue. In the ultra dense regime, other spectral quantities may be better candidates to look at than the maximal eigenvalue. This is no surprise: in \cite{touchette_equivalence_2015} it was shown that the relative entropy is the most sensitive global quantity to detect BEE, while other global quantities may detect BEE in certain settings and fail to do so in others. For instance, if the constraint is that the maximal eigenvalue takes a prescribed value, then clearly $\Delta_\infty = 0$ while BEE may still be possible.
\end{remark}


\paragraph{Outline.}

The remainder of this paper is organised as follows. In Section~\ref{sec:1} we recall the definition of the microcanonical and the canonical ensemble in the setting of constrained random graphs. Section~\ref{sec:method}  describes our main tool: a transfer method based on relative entropy, which carries over estimates on rare events from the canonical ensemble to the microcanonical ensemble, and describe its role in the general framework of BEE. In Section~\ref{sec:BEE} we prove Theorem~\ref{BEEspec}(1), in Section~\ref{sec:EE} we prove Theorem~\ref{BEEspec}(2).


\section{Gibbs ensembles for constrained random graphs}
\label{sec:1}

Consider the discrete probability space $\lep\GG_n, \mathcal{B}, \Pro\rip$, with $\GG_n$  the set of all simple graphs on $n$ vertices, $\mathcal{B}=2^{\GG_n}$ the power set of $\GG_n$ consisting of all the subsets of $\GG_n$, and $\Pro$  a probability measure. 

A \emph{constraint} is defined to be a vector-valued function $\vec{C}\colon\,\GG_n\to \Real^d$. Fix a value $\vec{C}^\star $ that is \emph{graphical}, i.e., $\vec{C}(g)=\vec{C}^\star $ for at least one $g\in\GG_n$. Define 
\begin{equation}
\Gamma_{\vec{C}^\star } = \left\{g\in \GG_n\colon\, \vec{C}(g)=\vec{C}^\star \right\}.
\end{equation}
The \emph{microcanonical ensemble} is the uniform probability distribution on $\Gamma_{\vec{C}^\star }$:
\begin{equation}
\Pmic(g) = \begin{cases}
1/|\Gamma_{\vec{C}^\star}|, 
&\text{if } g\in \Gamma_{\vec{C}^\star},\\
0, &\text{otherwise.}
\end{cases} 
\end{equation}
The \emph{canonical ensemble} is defined via the Hamiltonian $H(g,\vec{\theta})=\langle\vec{\theta},\vec{C}(g)\rangle$ (where $\langle\cdot,\cdot\rangle$ denotes the scalar product), namely,
\begin{equation}
\Pcan(g) = \frac{1}{Z_{\vec{\theta}^\star }}\,\ee^{-H(g,\vec{\theta}^\star)}, \qquad g \in \GG_n,
\end{equation} 
with the normalising factor $Z_{\vec{\theta}^\star } = \sum_{g\in \GG_n} \exp[-H(g,\vec{\theta}^\star)]$, called the \emph{partition function}. Note that both $\Pmic$ and $\Pcan$ depend on $n$, but we suppress this dependence.
The parameter $\vec{\theta}$ is set to the particular value $\vec{\theta}^\star$ that realises the constraint:
\begin{equation}
\label{eq:canconst}
\left. \Exp_{\mathrm{can}}\big[\vec{C}\big]\right|_{\theta=\theta^\star } = \vec{C}^\star .
\end{equation}
The constraint $\vec{C}^\star$, apart from being graphical, must also be \emph{irreducible}, i.e., no subset of the constraint is redundant \cite{squartini_garlaschelli_2020}. Once these conditions are met, the value of $\vec{\theta}^\star$ that solves \eqref{eq:canconst} is unique, and so the canonical ensemble is well defined (see the appendices in \cite{squartini_garlaschelli_2020} for further details).

The relative entropy of $\Pmic$ w.r.t.\ $\Pcan$ is defined as
\begin{equation}
\begin{split}
S_n(\Pmic\parallel\Pcan)
&= \sum_{g\in\GG_n}\Pmic(g) \log \frac{\Pmic(g)}{\Pcan(g)}
= \frac{1}{|\Gamma_{\vec{C}^\star }|} \sum_{g\in\Gamma_{\vec{C}^\star }} \log\frac{\Pmic(g)}{\Pcan(g)}\\
&= -\frac{1}{|\Gamma_{\vec{C}^\star }|}\log\big[|\Gamma_{\vec{C}^\star }|\,\Pcan(g^\star)\big]
\sum_{g\in\Gamma_{\vec{C}^\star }}1 = - \log\Pcan(\Gamma_{\vec{C}^\star })
\end{split}
\label{eq:relative}
\end{equation}
where we use the convention $0\log0=0$ and $g^\star$ is any graph in $\Gamma_{\vec{C}^\star}$. EE in the measure sense is defined as the vanishing of the relative entropy density, i.e., $\lim_{n\to\infty} n^{-1} S_n(\Pmic\parallel\Pcan)=0$ (see \cite{touchette_equivalence_2015}).

\section{Transfer method}
\label{sec:method}


\paragraph{Comparison of the two ensembles.} 

The additional freedom in the canonical ensemble implies that there is less dependence between the constituent random variables. In our case these random variables are the edges of the graph. For example, if the constraint is on the degree sequence, then the microcanonical ensemble corresponds to the \emph{hard configuration model} (which in the case of constant degrees becomes the regular random graph), while the canonical ensemble corresponds to the \emph{soft configuration model} (which is a special case of the generalized random graph model). The former requires an algorithm that randomly pairs half-edges and creates dependencies, while the latter is constructed via a sequence of independent random trials (which results in a multivariate Poisson-Binomial distribution for the degrees of the vertices~\cite{garlaschelli_covariance_2018}). Consequently, in the canonical ensemble calculations are carried out more easily. For example, a lot is known about spectral properties of adjacency matrices of random graphs under the canonical ensemble: because the entries of the adjacency matrix are independent, powerful tools from random matrix theory can be used. The challenge is to transfer properties from the canonical ensemble to the microcanonical ensemble without performing elaborate combinatorial computations.


\paragraph{Transfer principle.}
 
We start by noting that
\begin{equation}
\label{eq:quella}
\Pmic\lep B\rip 
= \frac{\Pcan(B)}{\Pcan(\Gamma_{\vec{C}^\star })}, 
\qquad B \subseteq \Gamma_{\vec{C}^\star }.
\end{equation}
The latter holds because $g \mapsto H(g,\vec{\theta}^\star )$ and $g \mapsto \Pcan(g)$ are constant on the support of $\Pmic$, i.e., all microcanonical realisations have the same probability under the canonical ensemble. In particular,
\begin{equation}
\Pcan(B\mid \Gamma_{\vec{C}^\star }) = \Pmic(B), \qquad B\in\mathcal{B},
\end{equation}
where again $\mathcal{B}=2^{\GG_n}$. Consequently, we have the following \emph{transfer principle}.
\begin{lemma}
For every $B\in\mathcal{B}$, if $\lim_{n\to\infty}\mathbb{P}_{\mathrm{can}}(B\mid \Gamma_{\vec{C}^\star })=0$, then $\lim_{n\to\infty} \mathbb{P}_{\mathrm{mic}}\lep B\rip=0$.
\end{lemma}


\paragraph{Distinguishing sets.}

Let $\mathcal{E}_\mathcal{P}\in\mathcal{B}$ be the subset of $\GG_n$ given by
\begin{equation}
\label{eq:event_p}
\mathcal{E}_\mathcal{P} = \left\{g\in\GG_n\colon\,g \text{ has property }\mathcal{P}\right\}.
\end{equation}
Write $[\mathcal{E}_\mathcal{P}]^c$ to denote the complementary event. The crucial step in the argument underlying the transfer method is to find the right event $[\mathcal{E}_\mathcal{P}]^c$ that asymptotically implies failure of the property $\mathcal{P}$ that we want to transfer from the canonical ensemble to the microcanonical ensemble.

For the remainder, two events are important: $\mathcal{E}_\mathcal{P}\cap\Gamma_{\vec{C}^\star}$ and $[\mathcal{E}_\mathcal{P}]^c\cap\Gamma_{\vec{C}^\star}$. These represent the sets that are in the support of $\Pmic$ for which property $\mathcal{P}$ holds and fails, respectively. Our focus will be on replacing $\Pcan([\mathcal{E}_\mathcal{P}]^c\cap\Gamma_{\vec{C}^\star})$ by $\Pcan([\mathcal{E}_\mathcal{P}]^c)$. Since $\Pmic([\mathcal{E}_\mathcal{P}]^c\cap\Gamma_{\vec{C}^\star}) \leq \Pmic([\mathcal{E}_\mathcal{P}]^c)$, if we are able to prove that $\lim_{n\to \infty} \Pmic([\mathcal{E}_\mathcal{P}]^c)=0$, then we also have $\lim_{n\to \infty} \Pmic([\mathcal{E}_\mathcal{P}]^c\cap\Gamma_{\vec{C}^\star})=0$, and we say that the property defining the set $\mathcal{E}_\mathcal{P}$ holds with high probability as $n\to\infty$. As explained in Section \ref{sec:1},
\begin{equation}
\label{eq:method}
\Pcan\big([\mathcal{E}_\mathcal{P}]^c \mid \Gamma_{\vec{C}^\star}\big) 
= \frac{\Pcan([\mathcal{E}_\mathcal{P}]^c\cap\Gamma_{\vec{C}^\star})}{\Pcan(\Gamma_{\vec{C}^\star})}
\leq \frac{\Pcan([\mathcal{E}_\mathcal{P}]^c)}{\Pcan(\Gamma_{\vec{C}^\star})},
\end{equation}
and so if we manage to prove that $\Pcan([\mathcal{E}_\mathcal{P}]^c)=o(\Pcan(\Gamma_{\vec{C}^\star}))$, then we obtain\\ $\lim_{n\to\infty} \Pmic([\mathcal{E}_\mathcal{P}]^c)=0$.


\paragraph{Role of relative entropy and BEE.}

Equation \eqref{eq:method} sets the scale at which the transfer method is effective. This scale is given by the denominator $\Pcan(\Gamma_{\vec{C}^\star })$. Indeed, if it happens that $\Pcan([\mathcal{E}_\mathcal{P}]^c)\neq o(\Pcan(\Gamma_{\vec{C}^\star }))$, then \eqref{eq:method} is ineffective. Importantly, from \eqref{eq:relative} we have
\begin{equation}
\label{eq:can}
\Pcan(\Gamma_{\vec{C}^\star }) = \ee^{-S_n(\Pmic\,\parallel\,\Pcan)}.
\end{equation}
This leads to an interesting connection between BEE and the transferability of a property $\mathcal{P}$: if $\Pcan([\mathcal{E}_\mathcal{P}]^c)=o(\ee^{-S_n(\Pmic\,\parallel\,\Pcan)})$, then $\lim_{n\to\infty} \Pmic([\mathcal{E}_\mathcal{P}]^c)=0$. Since EE coincides with $S_n(\Pmic\parallel\Pcan)=o(n)$, when the ensembles are equivalent it is easier to transfer. Our proof of Theorem~\ref{BEEspec}(2) makes use of precisely this fact, and $\mathcal{P}$ is a certain concentration inequality for the largest eigenvalue of the adjacency matrix. By contrast, BEE makes the transfer more difficult. Indeed, Theorem~\ref{BEEspec}(1) can be seen as an example where the same concentration inequality $\mathcal{P}$ cannot be transferred because the relative entropy is of higher order, namely, $S_n(\Pmic\parallel\Pcan)=\Theta(n\log n)$ \cite{squartini_breaking_2015,garlaschelli_covariance_2018}.


\paragraph{Largest eigenvalue.}

We know from the results in \cite{touchette_equivalence_2015} that ,whenever BEE occurs, there must exist quantities whose macrostate expectation is different under the two ensembles. Clearly, not all macroscopic quantities are good candidates for this. For instance, any linear combination of the constraints necessarily has the same expected value under the two ensembles. What we propose as a candidate is the largest eigenvalue of the adjacency matrix of the graph, because this is a highly nonlinear function of the imposed constraints and is sensitive to the global structure of the graph. In Sections~\ref{sec:BEE}--\ref{sec:EE} we will consider two examples of constraints in the dense regime: (1) fix the degrees of all the vertices; (2) fix the total number of edges. For the former we focus on the special case where all the degrees are equal. 

\begin{remark}
\label{rmk:generalization}
Since $\lambda_1(A) = \sup_{\norm{x}=1} x^T Ax$, Jensen's inequality implies that $\lambda_1(A)$ is a convex function of the entries of the matrix $A$, which means that for both ensembles 
\begin{equation}
\label{eq:jensen}
\lambda_1\lep\Exp_{(\cdot)}[A]\rip \leq \Exp_{(\cdot)}[\lambda_1(A)].
\end{equation}
Taking into account the results of Theorem \ref{thm:furedi} and Section \ref{sec:BEE}, we get
\begin{equation}
\label{eq:equalityjen}
\lambda_1(\Emic[A]) = \Emic[\lambda_1(A)] = \lambda_1(\Ecan[A]) \leq \Ecan[\lambda_1(A)].
\end{equation}
If, on top of the constraint on the degree sequence, we add more (compatible) constraints, then by exchangeability we still have $\lambda_1(\Emic[A])=\Emic[\lambda_1(A)]=\lambda_1(\Ecan[A])$. Applying \eqref{eq:jensen}, we therefore still expect that $\Emic[\lambda_1(A)] \leq \Ecan[\lambda_1(A)]$. This shows that $\lambda_1$ is particularly sensitive to the moments of the underlying degree sequence (as can also be seen from the \emph{power method} used in \cite{furedi_eigenvalues_1981,erdosSpectralStatisticsErdos2013}; see \eqref{eq:expansion} and \eqref{eq:lambda} below). We may therefore expect that our working hypothesis holds in all those cases where BEE forces the degree sequence to assume either a different mean of a different variance in the two ensembles, as in the case under study.
\end{remark}

\section{Proof of Theorem~\ref{BEEspec}(1): constraint on the degree sequence}
\label{sec:BEE}

In what follows we suppress the $n$-dependence from $p(n),d(n),\lambda_1(n)$, writing $p,d,\lambda_1$. The $d$-regular random graph with $n$ vertices, written $G_{n,d}$, coincides with the microcanonical ensemble with constraint $\vec{C}^\star =(d,\ldots,d)$ on the degree sequence, where we allow $d=d(n)$. The largest eigenvalue of the adjacency matrix of $G_{n,d}$ equals $d$, irrespective of $n$. The \ER random graph with retention probability $p=d/(n-1)$ coincides with the canonical ensemble with the same constraint.
 
In order to understand the difference in behaviour of $\lambda_1$ under the two ensembles, we need Theorem \ref{thm:furedi}. Indeed, the result in \eqref{CLT}, which actually holds for a generic symmetric random matrix subject to specific regularity conditions, can be interpreted as follows. The adjacency matrix $A$ associated with $G(n,p)$ consists of elements $\{a_{ij}\}_{i,j\in[n]}$ that are identically $0$ when $i=j$ and \emph{Bernoulli random trials} ($a_{ij}=0,1$) with success probability $p$ when $i\neq j$. The largest eigenvalue of the deterministic matrix $\bar{A}$ whose entries are $\bar{a}_{ij}=\Exp_{\mathrm{can}}[a_{ij}]=p$ when $i\neq j$ and $\bar{a}_{ij}=0$ when $i=j$ is given by $\lambda_1(\bar{A})= (n-1)p$. Hence, compared to $\lambda_1(\bar{A})$, $\lambda_1$ is shifted by a random variable whose expected value is $(1-p)$ and is distributed as $\mathcal{N}\lep 1-p,2p(1-p)\rip$ under certain conditions on $d$ (see \cite[equation 6.10]{erdosSpectralStatisticsErdos2013}) plus an error term of order dependent on the considered regime ( $O(1/\sqrt{n})$ for $p$ constant). It is important to note that the parameters of this shift depend on $p$ only. In \cite{furedi_eigenvalues_1981,erdosSpectralStatisticsErdos2013} it is shown that \eqref{CLT} relies on the fact that in the canonical ensemble the eigenvector $\vec{v}_1$ corresponding to the largest eigenvalue $\lambda_1$ is very close to the vector $\ones=(1,\ldots,1)$ (i.e., the norm of the projection of $\vec{v}_1$ onto $\ones$ is much larger than the norm of the projection of $\vec{v}_1$ onto the perpendicular space $\ones^\perp$). 

It was shown in \cite{garlaschelli_covariance_2018} that BEE holds in the all regimes covered in Theorem~\ref{BEEspec}(1), namely the delta tame regime, which corresponds to $\delta \leq p = d/(n-1) \leq 1-\delta$ with $\delta\in(0,\frac{1}{2}]$ (see \cite[Definition 1.1]{garlaschelli_covariance_2018}) and the sparse regime ($d=o(\sqrt{n})$). Hence the claim in Theorem~\ref{BEEspec}(1) follows.


\section{Proof of Theorem~\ref{BEEspec}(2): constraint on the total number of edges}
\label{sec:EE}

Consider the case where the constraint is on the total number of edges: $\vec C(g) = \vec C^\star  = \binom{n}{2}p$ for some $p \in (0,1)$. Then the canonical ensemble is still the \ER random graph with parameter $p$. It was proved in \cite{squartini_breaking_2015} that the two ensembles are asymptotically equivalent on scale $n$. In particular, it was shown that $S_n(\Pmic\parallel\Pcan) = \log n + \Theta(1)$. The canonical probability of drawing a microcanonical realization is given by \eqref{eq:can}:
\begin{equation}
\label{relentbd}
\Pcan(\Gamma_{\vec{C}^\star }) = \ee^{-S_n(\Pmic\,\parallel\,\Pcan)} =\ee^{-\log n+\Theta(1)} = \Theta(n^{-1}).
\end{equation}
Together with \eqref{eq:method}, this tells us that if we can find an event $[\mathcal{E}_\mathcal{P}]^c$ such that $\Pcan([\mathcal{E}_\mathcal{P}]^c)=o(n^{-1})$, then we know that $\lim_{n\to\infty} \Pmic(\mathcal{E}_\mathcal{P})=1$. Our goal is to use the results in \cite{erdosSpectralStatisticsErdos2013} to apply \eqref{eq:method} with \eqref{relentbd}.

In Section~\ref{sec:erdos} we show how our results follow from \cite[Theorem 6.2]{erdosSpectralStatisticsErdos2013} both in the dense and the non-dense regime. In Section~\ref{sec:ced} and \ref{sec:cele} we focus on the dense regime and show how our results follow by making the concentration inequalities used in \cite{furedi_eigenvalues_1981} tighter. In particular, we will find that the approach heavily depends on the ability of identifying good concentration inequalities for the degree sequence, which is a special case of the bounds presented in \cite{erdosSpectralStatisticsErdos2013}. The heavy dependence on the degree sequence is further evidence of what was said in Remark~\ref{rmk:generalization}. In Section~\ref{sec:ced} we prove a concentration inequality for the degrees under the canonical ensemble (Lemma~\ref{lem:conc}) that is of independent interest. In Section~\ref{sec:cele} we use this to prove a concentration inequality for a functional of the degrees that approximates the largest eigenvalue well in the dense regime (Lemma~\ref{lem:mean}). In Section~\ref{sec:tr} we transfer the results from the previous sections to the microcanonical ensemble (Lemma~\ref{lem:tr}), and show that this leads to a negligible shift of the expected largest eigenvalue.

\subsection{Proof of Theorem~\ref{BEEspec}(2) via \cite{erdosSpectralStatisticsErdos2013}}
\label{sec:erdos}

In \cite[Chapter 6]{erdosSpectralStatisticsErdos2013} the largest eigenvalue of matrices of the form $A=A_0+\Exp[ A]=A_0+f\ket{\vec{e}\,}\bra{\vec{e}\,}$ is studied, where $A_0$ is a matrix with mean-zero entries, $\vec{e}=\frac{1}{\sqrt{n}}\lep1,1, \dots,1\rip^T$ and $1+\varepsilon_0\leq f \leq N^C $, with $\varepsilon_0,C\in(0,\infty)$ constants. In order to consider only adjacency matrices $A_s$ of simple graphs, we have to get rid of the diagonal of $A$. This can be easily done by considering $A_s=A-pI= A_0+f\ket{\vec{e}\,}\bra{\vec{e}\,}-pI$, where $p=p(n)$ is the retention probability that appears Theorem~\ref{thm:furedi} subject to \eqref{regime} and $f=np$. We note that if $\lambda_1$ is the largest eigenvalue of $A$, then $\lambda_1-p$ is the largest eigenvalue of $A_s$, so it suffices to study the largest eigenvalue of $A$.

Let $\tilde{A}$ be the normalized version of $A$, defined by $\tilde{A}=A/\sqrt{np(1-p)}$. This scaling is needed in order to have $\norm{H}=O(1)$ with high probability. Let $\tilde{A}_0$ be the centered version of $\tilde{A}$, i.e., $\tilde{A}_0=\tilde{A}- \Ecan\left[{\tilde{A}}\right]$. It is easy to see that $\Ecan\left[{\tilde{A}}\right]$ can be expressed as $\sqrt{np/(1-p)}\ket{\vec{e}\,}\bra{\vec{e}\,}$, where again $\vec{e}=\frac{1}{\sqrt{n}}\lep1,1, \dots,1\rip^T$. Following \cite[Theorem 6.2]{erdosSpectralStatisticsErdos2013}, we say that an event $\mathcal{E}$ holds with $(\xi,\nu)$-high probability when
\begin{equation}
\label{eq:nuxihigh}
\Pro(\mathcal{E}^c)\leq \ee^{-\nu (\log n)^\xi},
\end{equation}
where $\nu$ and $\xi$ can be two positive $n$-dependent constants bounded from below by $\nu>0$ and $\xi>1$. Note that $e^{-\nu (\log n)^\xi }=o(n^{-1})$ whenever $\nu>0$ and $\xi>1$. Thus, if an event $\mathcal{E}_\mathcal{P}$ of the type described in \eqref{eq:event_p} holds with $(\xi,\nu)$-high probability under $\Pcan$, then by \eqref{eq:method} and \eqref{relentbd} $\mathcal{E}_\mathcal{P}$ it holds also under $\Pmic$. Starting from the equation 
\begin{equation}
\lep\mathbb{I}-\frac{\tilde{A}_0}{\lambda_1}\rip\lambda_1 \vec{v}=\sqrt{\frac{np}{1-p}}\langle \vec{e}, \vec{v} \rangle\vec{e},
\end{equation}
where $ \vec{v} $ is the eigenvector associated with $\lambda_1$ and $\mathbb{I}$ is the identity matrix, after multiplying by $(\mathbb{I}-\frac{\tilde{A}_0}{\lambda_1})^{-1}$ and projecting on $\vec{e}$ we obtain the following series for $\lambda_1$:
\begin{equation}
\label{eq:expansion}
\lambda_1 = \sqrt{\frac{np}{1-p}}
\sum_{k\in\mathbb{N}_0} \left\langle \vec{e}, \lep\frac{\tilde{A}_0}{\lambda_1}\rip^k\vec{e}\right\rangle. 
\end{equation}
We see that for the series to converge we need $\norm{\tilde{A}_0}/\lambda_1<1$. From \cite[Lemma 4.3]{erdosSpectralStatisticsErdos2013} (see also \cite{alon_concentration_2002,soshnikov_universality_1999,tao_topics_2012,vu_spectral_2007}) and the leading order of \eqref{eq:expansion} (see also \cite[Eq.(6.5)]{erdosSpectralStatisticsErdos2013}) we have that $\norm{\tilde{A}_0}/\lambda_1<1$ with $(\xi,\nu)$-high probability (which also holds for the microcanonical ensemble). Iterating \eqref{eq:expansion}, we get that with $(\xi,\nu)$-high probability
\begin{equation}
\label{eq:expansion2}
\begin{aligned}
&\lambda_1 = \sqrt{\tfrac{np}{1-p}}+\left\langle \vec{e}, \tilde{A}_0
\vec{e}\right\rangle+\lep\left\langle \vec{e}, \tilde{A}^2_0\vec{e}\right\rangle-\left\langle \vec{e}, \tilde{A}_0\vec{e}\right\rangle^2\rip\lep\sqrt{\tfrac{np}{1-p}}\rip^{-1} \\
&+\lep\left\langle \vec{e}, \tilde{A}_0\vec{e}\right\rangle^3-3\left\langle \vec{e}, \tilde{A}_0\vec{e}\right\rangle\left\langle \vec{e}, \tilde{A}^2_0\vec{e}\right\rangle\rip\lep\sqrt{\tfrac{np}{1-p}}\rip^{-2}
+ O\lep\lep\sqrt{\tfrac{np}{1-p}}\rip^{-3}+\lep\tfrac{(np)q}{1-p}\rip^{-1}\rip,
\end{aligned}
\end{equation}
where $q$ is the parameter defined in Theorem \ref{thm:furedi}. Taking expectations, using \cite[Lemma 6.5]{erdosSpectralStatisticsErdos2013} and scaling back, we get \eqref{CLT} for $\Exp[\lambda_1]-p$, the expected eigenvalue of $A_s$. Note that all the bounds hold with $(\xi,\nu)$-high probability. We can therefore conclude via \eqref{eq:method} that \eqref{CLT} approximates $\lambda_1$ with a vanishing error also in the microcanonical ensemble, where the constraint is on the total number of edges. Together with the result of Lemma~\ref{lem:mean}, we conclude that $\lim_{n\to\infty} \lep\Exp_{\mathrm{can}}\left[\lambda_1\right] - \Exp_{\mathrm{mic}}\left[\lambda_1\right]\rip = 0$. 


\subsection{Concentration for the degrees under the dense canonical ensemble}
\label{sec:ced}

For the remainder of the paper we take $p\in(0,1)$ constant and $A$ to be the unnormalized adjacency matrix. For $i\neq j$, $\Exp_{\mathrm{can}}[a_{ij}]=p$ and $\Var_\mathrm{can}[a_{ij}]= p(1-p)$. In what follows we abbreviate $\mu=p$ and $\sigma^2= p(1-p)$. We write $\ones= \vec{v}_1 + \vec{r}$ with $\vec{r}\in\ones^\perp$, $\langle\vec{v}_1, \vec{r}\rangle=0$ and $A\vec{v}_1=\lambda_1\vec{v}_1$. Following the power method in \cite{meyer_matrix_2000}, we define 
\begin{equation}
\vec{K} = A\ones=A(\vec{v}_1 + \vec{r})=\lambda_1 \vec{v}_1 + A\vec{r},
\end{equation}
which is the vector of row sums of the matrix $A$, i.e., the vector of degrees of the vertices (the degree sequence). Centering $\vec{K}$ by $\Theta \ones$ with $\Theta=\Exp[K_i]=(n-1)p$ and using $\ones= \vec{v}_1 + \vec{r}$, we get
\begin{equation}
\vec{K}-\Theta\ones= (\lambda_1-\Theta) \vec{v}_1 + (A\vec{r}-\Theta\vec{r}).
\end{equation}

Our key step is the following lemma.

\begin{lemma}
\label{lem:conc}
With $\sigma^2$ denoting $p(1-p)$, there exist two constants $c_1,c_2\in (0,\infty)$ such that
\begin{equation}
\label{eq:conc}
\mathbb{P}_{\mathrm{can}}\lep\abs{\sum_{i=1}^n \lep K_i-\Theta \rip^2 -\sigma^2n(n-1)} \geq t\rip 
\leq c_2 \ee^{-c_1t/n^{3/2}}.
\end{equation}
\end{lemma}

\begin{proof}
The term $\sum_{i=1}^n \lep K_i-\Theta \rip^2$ can be written as 
\begin{equation}
\label{eq:sum}
\sum_{i=1}^n \lep \sum_{j=1}^n \lep a_{ij}-\Exp_{\mathrm{can}}[a_{ij}]\rip \rip^2
=\sum_{i=1}^n \lep \sum_{j=1}^n  b_{ij}\rip^2
=\sum_{i=1}^n  \sum_{j=1}^n\sum_{k=1}^n b_{ij}b_{ik},
\end{equation}
where
\begin{equation}
b_{ij} =a_{ij}-\Exp_{\mathrm{can}}[a_{ij}]=\begin{cases}
a_{ij}-p, &\quad \text{if } i\neq j,\\
0, &\quad \text{if }i=j, 
\end{cases}
\end{equation}
are the centred entries of the adjacency matrix. Note that 
\begin{equation}
\Exp_{\mathrm{can}}\left[\sum_{i=1}^n  \sum_{j=1}^n\sum_{k=1}^n b_{ij}b_{ik}\right] = \sigma^2n(n-1).
\end{equation} 

Straightforward counting shows that the sum in \eqref{eq:sum} contains $O(n^3)$ different terms. Let us represent $b_{ij}=b_{ji}$ by a variable $X_\alpha$, $\alpha\in\left[\binom{n}{2}\right]$. Then \eqref{eq:sum} can be rewritten in the form
\begin{equation}
\label{eq:sum2}
\sum_{\alpha, \beta \in \left[\binom{n}{2}\right]} h_{\alpha\beta}X_\alpha X_\beta,
\end{equation}
which is the quadratic form of the matrix $H=\{h_{\alpha\beta}\}_{\alpha, \beta \in \left[\binom{n}{2}\right]}$. Because there is a one-to-one correspondence between the terms in \eqref{eq:sum2} and \eqref{eq:sum}, we can conclude that $H$ has $O(n^3)$ entries, whose values are either $1$ (off-diagonal) or $2$ (diagonal). We can apply to \eqref{eq:sum2} the \emph{Hanson-Wright} inequality (see \cite{hanson_bound_1971} or \cite[Theorem 1.4, item 6]{adamczak_concentration_2015}).

\begin{theorem}
\label{thm:ineq}
Let $X=(X_1,\ldots X_N)$ be mean-zero square-integrable random variables taking values in $\Real$, and let $\xi>0$ be such that
\begin{equation}\label{condxi}
\norm{X}_{\psi_2}=\inf\left\{t>0 : \Exp\left[\exp\lep\|X\|_2^2/t^2\rip\right]\leq 2\right\}\leq \xi.
\end{equation}
Let $H=(h_{\alpha\beta})_{\alpha\beta\in [N]}$ be a real symmetric matrix. Then $Y=\sum_{\alpha,\beta \in [N]} h_{\alpha\beta} X_\alpha X_\beta$ satisfies
\begin{equation}
\label{normmin}
\Pro\lep\abs{Y-\Exp[Y]}\geq t \rip
\leq 2\exp\lep-\frac{1}{C}\min\left\{\frac{t^2}{\xi^4\norm{H}^2_{\mathrm{HS}}},
\frac{t}{\xi^2 \norm{H}_{\ell_2^N\to\ell_2^N}}\right\}\rip, \qquad t>0,
\end{equation}
where $C$ is a suitable constant, $\norm{H}^2_{\mathrm{HS}}=\sum_{\alpha,\beta\in[N]} h_{\alpha\beta}^2$ is the Hilbert-Schmidt norm of $H$, and 
\begin{equation}
\norm{H}^2_{\ell_2^N\to\ell_2^N}=\sup\left\{\sum_{\alpha,\beta\in[N]}h_{\alpha\beta}x_\alpha y_\beta: \sum_{\alpha\in [N]} x_\alpha^2\leq1, \sum_{\alpha\in [N]} y_\alpha^2\leq1\right\}
\end{equation}
is the $\ell^N_2 \to \ell^N_2$ norm of $H$.
\end{theorem}

In our setting, $N=\binom{n}{2}$. Since $\abs{X_\alpha}<1$, we have $\norm{X}_{\psi_2} \leq 1/\log 2$, so that \eqref{condxi} applies with $\xi=1/\log 2$. Since $H$ has bounded entries, we have $\norm{H}^2_{\mathrm{HS}}=O(n^3)$. Moreover, by the Cauchy-Schwarz inequality we have
\begin{equation}
\norm{H}^2_{\ell_2^N\to\ell_2^N} = \sup\left\{\norm{Hx}_2\colon\,\norm{x}_2\leq1\right\} = \norm{H}_{\mathrm{op}},
\end{equation} 
where the latter is the operator norm of $H$. But
\begin{equation}
\norm{H}_{\mathrm{HS}}^2 = \text{Tr}(H^\dagger H)\geq\lambda_\text{max}(H^\dagger H)=\norm{H}^2_{\mathrm{op}},
\end{equation}
and so the exponent in the right-hand side of \eqref{normmin} is bounded below by  
\begin{equation}
\min\left\{\frac{t^2}{\xi^4n^3},\frac{t}{\xi^2n^{3/2}}\right\} \geq \frac{c_3 t}{n^{3/2}},
\end{equation}
where $c_3$ is a suitable constant. Taking $c_1\leq{c_3}/{C}$, with $C$ the constant appearing in  \eqref{normmin}, we obtain \eqref{eq:conc}. 
\end{proof}

We end this section with an immediate consequence of Lemma \ref{lem:conc}. Picking $t=\sigma^2n^2$ and using that, for appropriately chosen constants $C_1, C_2, C_3,C_4$,
\begin{equation}
\frac{\sigma^4n^4}{\norm{H}^2_{\mathrm{HS}}} \geq \frac{\sigma^4 n^4}{C_1 n^3} \geq C_2 n, \qquad
\frac{\sigma^2 n^2}{L^2 \norm{H}_{\mathrm{op}}}\geq \frac{\sigma^2n^2}{C_3 \norm{H}_{\mathrm{HS}}}
\geq C_4\sqrt{n},
\end{equation}
we find that there are constants $\tilde c\leq {C_4}/{C}$ and $\tilde C$ such that
\begin{equation}
\label{eq:conc2}
\begin{aligned}
&\Pcan\lep \abs{\sum_{i=1}^n \lep K_i-\Theta \rip^2 -\sigma^2n^2} \geq 2\sigma^2n^2\rip\\ 
&\qquad \leq 2\exp\lep-\frac{1}{C}\min\left\{\frac{4 \sigma^4n^4 }{\norm{H}^2_{HS}},
\frac{2\sigma^2n^2}{\norm{H}_{op}}\right\}\rip \leq \tilde{C} \ee^{-\tilde{c}\sqrt{n}}.
\end{aligned}
\end{equation}


\subsection{Concentration for the largest eigenvalue under the dense canonical ensemble}
\label{sec:cele}

After applying $A$ once to $\ones$, we must find a suitable normalization in order to isolate $\lambda_1$. This is given by 
\begin{equation}
\label{eq:boh}
\frac{\sum_{i=1}^n K_i^2}{\sum_{i=1}^n K_i} 
= \frac{\langle \vec{K},\vec{K} \rangle}{\langle \ones, \vec{K}\rangle}
= \frac{\|A\ones\|}{\langle\ones, A \ones\rangle}
= \lambda_1+ \frac{\|A\vec{r}\|^2-\lambda_1\langle\vec{r},A\vec{r}\rangle}{\sum_{i=1}^n K_i}.
\end{equation}

In \cite{furedi_eigenvalues_1981}, it was shown that $\sum_{i=1}^n K_i^2/\sum_{i=1}^n K_i$ approximates $\lambda_1$ with high probability, in the sense that for any $x>0$,
\begin{equation}
\label{eq:concent}
\Pcan\lep\abs{\frac{\sum_{i=1}^n K_i^2}{\sum_{i=1}^n K_i}-\frac{\sum_{i=1}^n K_i}{n}-\frac{\sigma^2}{\mu}} 
\geq \frac{3\sigma^2x}{\sqrt{n}}\rip \leq \frac{1}{x^2}, 
\end{equation}
which with the choice $x=\sqrt{n}$ leads to an upper bound of order $1/n.$ As it turns out, however, in order to transfer the estimates   to the microcanonical ensemble via \eqref{eq:method}, we need the upper bound to hold with probability $o(1/n)$. This result is covered by  the following lemma.

\begin{lemma}
\label{lem:mean}
Let $\vec{K}$ be as before, and $\mu=p$, $\sigma^2=p(1-p)$. For any $\gamma>0$ there exist $\gamma', \gamma_1, \gamma_2$ satisfying $c_1\gamma_1, \gamma_2>1$, with $c_1$ the constant in \eqref{eq:conc}, such that
\begin{equation}
\label{eq:inequality}
\mathbb{P}_{\mathrm{can}}\lep\abs{\frac{\sum_{i=1}^n K_i^2}
{\sum_{i=1}^n K_i}-\frac{\sum_{i=1}^n K_i}{n}-\frac{\sigma^2}{\mu}} 
\geq \frac{\gamma}{\sqrt{n}}\rip
\leq \frac{\gamma'}{n^{\min\{c_1\gamma_1,\gamma_2\}}}.
\end{equation}
\end{lemma}

\begin{proof}
First note that
\begin{equation}
\Exp_{\mathrm{can}}\left[\frac{\sum_{i=1}^n K_i}{n}\right]
=\frac{1}{n}\sum_{i=1}^n \Exp_{\mathrm{can}}[K_i]
= (n-1)p = \Theta
\end{equation}
and write 
\begin{equation}
\label{quotients}
\frac{\sum_{i=1}^n K_i^2}{\sum_{i=1}^n K_i}-\frac{\sum_{i=1}^n K_i}{n}
= \frac{\sum_{i=1}^n \lep K_i-\Theta \rip^2}{\sum_{i=1}^n K_i}
-\frac{\lep n^{-1}\sum_{i=1}^n K_i-\Theta\rip^2}{n^{-1}\sum_{i=1} K_i}.
\end{equation}

To analyse the first ratio in \eqref{quotients}, note that
\begin{equation}
\sum_{i=1}^n K_i 
=\sum_{i,j\in[n]} a_{ij}
= 2\sum_{i,j\in[n], j>i} a_{ij}.
\end{equation}
Applying Hoeffding's inequality (see e.g.~\cite{bercu_concentration_2015,boucheron_concentration_2013}), we have 
\begin{equation}
\label{eq:hoeff}
\Pcan\lep\abs{\sum_{i,j\in[n], j>i} a_{ij}-\frac{n(n-1)}{2}\mu}\geq t\rip \leq 2\exp\lep-\frac{4t^2}{n(n-1)}\rip.
\end{equation}
Take $t=n\sqrt{\gamma_2\log n}$ in \eqref{eq:hoeff} with $\gamma_2>1$ and apply Lemma~\ref{lem:conc} with $t=n^{3/2} \gamma_1\log n$, with $\gamma_1c_1>1$ and $c_1$ the constant in the exponential bound of \eqref{eq:conc}. Then, for some $\gamma>0$,
\begin{equation}
\frac{\sum_{i=1}^n \lep K_i-\Theta \rip^2}{\sum_{i=1}^n K_i}
\leq\frac{n(n-1)\sigma^2+n^{3/2}\gamma_1\log n}{n(n-1)\mu+n\sqrt{\gamma_2\log n}}
\leq\frac{\sigma^2}{\mu} + \frac{\gamma}{\sqrt{n}}
\end{equation}
with probability at least $1-{1}/{n^{\gamma_1c_1}}-{1}/{n^{\gamma_2}}$. Similarly, the probability of 
\begin{equation}
\frac{\sum_{i=1}^n \lep K_i-\Theta \rip^2}{\sum_{i=1}^n K_i}
\geq \frac{\sigma^2}{\mu} - \frac{\gamma}{\sqrt{n}}
\end{equation}
is bounded from below by $1-{1}/{n^{\gamma_1c_1}}-{1}/{n^{\gamma_2}}$. Hence
\begin{equation}
\Pcan\lep\abs{ \frac{\sum_{i=1}^n \lep K_i-\Theta \rip^2}{\sum_{i=1}^n K_i}-\frac{\sigma^2}{\mu}}
\geq \frac{\gamma}{\sqrt{n}} \rip\leq \frac{\gamma'}{n^{\min\{\gamma_1c_1,\gamma_2\}}}.
\end{equation}
To analyse the second ratio in \eqref{quotients}, we write 
\begin{equation}
\lep n^{-1} \sum_{i=1}^n K_i-\Theta\rip^2 
=\frac{1}{n^2}\lep 2\sum_{i,j\in[n], j>i} \lep a_{ij}
-\Exp_{\mathrm{can}}[a_{ij}]\rip\rip^2,
\end{equation}
and apply Hoeffding's inequality with $t=O(n^2)$ twice. This gives
\begin{equation}
\Pcan\lep\frac{\lep n^{-1}\sum_{i=1}^n K_i-\Theta\rip^2}{n^{-1}\sum_{i=1}^n K_i} > \frac{\tilde\gamma}{n}\rip
\leq \tilde\gamma_2e^{-\tilde\gamma_1 n^2}, \qquad \tilde\gamma>0,
\end{equation}
where $\tilde\gamma_2$ and $\tilde\gamma_1$ are suitable constants. Applying the union bound to the complementary events, we obtain \eqref{eq:inequality}.
\end{proof}


\subsection{Transfer to the dense microcanonical ensemble}
\label{sec:tr}

Next we use the transfer method to pass the property characterised by the event in \eqref{eq:inequality} to the microcanonical ensemble. Indeed, using the notation of Section~\ref{sec:method}, we identify
\begin{equation}
\abs{\frac{\sum_{i=1}^n K_i^2}{\sum_{i=1}^n K_i}-\frac{\sum_{i=1}^n K_i}{n}-\frac{\sigma^2}{\mu}}
\geq \frac{\gamma}{\sqrt{n}}
\end{equation}
as the event $\mathcal{E}_\mathcal{P}^c$, i.e., the set of graphs that do not possess the property that we would like to pass on. The fact that $\Pcan\lep \mathcal{E}_\mathcal{P}^c\rip$ tends to zero faster than $\Pcan(\Gamma_{\vec{C}^\star })$ (as $n\to\infty$, that is) tells us that also $\Pmic\lep \mathcal{E}_\mathcal{P}^c\rip$ tends to zero, and implies that
\begin{equation}
\lim_{n\to\infty} \Pmic\lep\abs{\frac{\sum_{i=1}^n K_i^2}{\sum_{i=1}^n K_i}-\frac{\sum_{i=1}^n K_i}{n}-\frac{\sigma^2}{\mu}} \leq \frac{\gamma}{\sqrt{n}}\rip = 1.
\end{equation} 

Thus, in the microcanonical ensemble $\sum_{i=1}^n K_i^2/\sum_{i=1}^n K_i$ concentrates around the sum $n^{-1}\sum_{i=1}^n K_i+{\sigma^2}/{\mu}$ with an error of order $1/\sqrt{n}$. However, we need to also see what $n^{-1}\sum_{i=1}^n K_i+{\sigma^2}/{\mu}$ is in the microcanonical ensemble. The term ${\sigma^2}/{\mu}$, a constant equal to $1-p$, is in accordance with the constraint in the microcanonical ensemble. For the other term we have $n^{-1} \sum_{i=1}^n K_i =(n-1)p$. The two together give precisely the expected value in the canonical ensemble, as follows from Proposition~\ref{thm:furedi}. Hence we only need to show that $\sum_{i=1}^n K_i^2/\sum_{i=1}^n K_i$ concentrates around $\lambda_1$ also in the microcanonical ensemble, for which we can once more use the transfer method.

\begin{lemma}
\label{lem:tr}
For any $\eta>0$, there exist $\zeta$ and $\Lambda$ such that
\begin{equation}
\label{eq:lambda}
\mathbb{P}_{\mathrm{can}}\lep\abs{\frac{\sum_{i=1}^n K_i^2}{\sum_{i=1}^n K_i}-\lambda_1} 
\geq \frac{\eta}{\sqrt{n}}\rip \leq\Lambda e^{-\zeta \sqrt{n}} .
\end{equation}
\end{lemma}

\begin{proof}
We need to show that the last term in \eqref{eq:boh},
\begin{equation}\label{eq:minimizing}
\frac{\norm{A\vec{r}}^2-\lambda_1\langle\vec{r},A\vec{r}\rangle}{\sum_{i=1}^n K_i},
\end{equation}
is small. First we show that $\norm{\vec{r}}$ is bounded in probability. Indeed,
\begin{equation}
\sum_{i=1}^n \lep K_i-\Theta \rip^2=(\lambda_1-\Theta)^2\norm{\vec{v_1}}^2+ \norm{A\vec{r}-\Theta \vec{r}}^2.
\end{equation}
Since $(\lambda_1-\Theta)^2\norm{\vec{v_1}}^2\geq0$, we have $\norm{A\vec{r}-\Theta \vec{r}}^2<\sum_{i=1}^n \lep K_i-\Theta \rip^2$. By the Courant-Fisher theorem \cite[Theorem 1.3.2]{tao_topics_2012}, we get that $\norm{A\vec{r}-\Theta \vec{r}} \geq \abs{\Theta - \lambda_2}\norm{\vec{r}}$ (indeed, $\abs{\Theta - \lambda_i}\geq\abs{\Theta - \lambda_2}$ for $i>2$). Next, we need a concentration inequality for $\lambda_2$. Use \cite[Theorem 1]{alon_concentration_2002} plus the fact that the largest eigenvalue of a centred matrix is of order $O(\sigma \sqrt{n})$ almost surely \cite{soshnikov_universality_1999}, \cite[Theorem 2.3.24]{tao_topics_2012}, \cite[Theorem 1.3]{vu_spectral_2007}. Again use the Courant-Fisher theorem to pass to the non-centred case \cite[Lemma 1]{furedi_eigenvalues_1981}. We find that, for any $\beta>2$ and for $\tilde n$ large enough,
\begin{equation}
\label{eq:second}
\Pcan\lep \max_{i>1} \abs{\lambda_i} \geq \beta\sigma\sqrt{n}\rip \leq 4e^{-\zeta_1n},\quad n>\tilde n
\end{equation}
where $\zeta_1$ is a suitable constant. Since $\max_{i>1} \abs{\lambda_i}\geq \lambda_2\geq0$, we can bound $\lambda_2\leq \beta\sigma\sqrt{n}$ with high probability. Using \eqref{eq:conc2}, we have
\begin{equation}
\Pcan\left(\norm{\vec{r}}^2<\frac{\lep K_i-\Theta \rip^2}{(\Theta-\lambda_2)^2}
<\frac{n^2\sigma^2}{(\mu n-\beta\sigma\sqrt{n})^2}<\frac{4\sigma^2}{\mu^2}\right)
\geq 1-4e^{-\zeta_1n} - \tilde{C}e^{-\tilde{c}\sqrt{n}},
\end{equation}
as a consequence of the union bound applied to the last term of $\Pro \lep\cap_n {\mathcal E}_n \rip=1-\Pro(\cup_n [{\mathcal E}_n]^c)$, with $[{\mathcal E}_n]^c$ denoting the events described by Lemma \ref{lem:conc} and \eqref{eq:second}. Thus, we have 
\begin{equation}
\label{eq:rnorm}
\Pcan\lep\norm{\vec{r}}^2 \geq \frac{4\sigma^2}{\mu^2}\rip \leq \tilde C_1 e^{-\tilde c_1 \sqrt{n}},
\end{equation}
where $\tilde C_1 $ and $\tilde c_1$ are suitable constants.

All the other terms in \eqref{eq:minimizing} can be obtained by repeatedly using \eqref{eq:rnorm}, \eqref{eq:second} and \eqref{eq:conc2}. Note that in order to get \eqref{eq:rnorm} we have used both \eqref{eq:second} and \eqref{eq:conc2}, and the events that these inequalities identify.  Thus, using \eqref{eq:second} twice, we obtain
\begin{equation}
\Pcan\left(\norm{A\vec{r}}^2\leq\lambda_2^2 \norm{\vec{r}}^2\leq \frac{50\,\sigma^4}{\mu^2}n\right)\geq 
1-4e^{-\zeta_1n} - \tilde{C}e^{-\tilde{c}\sqrt{n}}.
\end{equation}
Therefore  
\begin{equation}
\Pcan\lep\norm{A\vec{r}}^2 \geq \frac{50\,\sigma^4}{\mu^2}n\rip  \leq \tilde C_2 e^{-\tilde c_2 \sqrt{n}},
\end{equation}
where $\tilde C_2 $ and $\tilde c_2$ are suitable constants. In the same way we can bound $\abs{\langle\vec{r},A\vec{r}\rangle}\leq\norm{\vec{r}}\norm{A\vec{r}}$, which yields 
\begin{equation}
\Pcan\lep\abs{\langle\vec{r},A\vec{r}\rangle} \geq \frac{2\sqrt{50}\sigma^3}{\mu^2}\sqrt{n}\rip
\leq \tilde C_3 e^{-\tilde c_3 \sqrt{n}}.
\end{equation}
Now, using the trivial deterministic bound $\lambda_1\leq \max_i\sum_j\abs{a_{ij}}<n$ and Hoeffding's inequality on $\sum_{i=1}^n K_i =2\sum_{j>i} a_{ij}$, we can conclude that, for any $\eta>0$, 
\begin{equation}
\Pcan\lep\abs{ \frac{\norm{A\vec{r}}^2-\lambda_1\langle\vec{r},A\vec{r}\rangle}{\sum_{i=1}^n K_i}}
\geq \frac{\eta}{\sqrt{n}}\rip \leq  \Lambda e^{-\zeta \sqrt{n}},
\end{equation}
where $\zeta$ and $\Lambda$ are suitable constants. Thus, recalling \eqref{eq:boh}, we have settled \eqref{eq:lambda}.
\end{proof}

We thus find that the probability in the canonical ensemble of the event in \eqref{eq:lambda} is $o\lep1/n\rip$, which confirms the results of Section~\ref{sec:erdos}. In particular, we have shown that the central object is the ratio  $\sum_{i=1}^n K_i^2/\sum_{i=1}^n K_i$.

\begin{remark}
{\rm The constants in the right-hand side of \eqref{eq:inequality} can be chosen freely. By Lemma \ref{lem:tr}, this means that for any choice of constraint for which $S_n(\Pmic\parallel\Pcan)=O(\log n)$ and the canonical ensemble is the \ER random graph, we have that $\lambda_1$ is close to $n^{-1}\sum_{i=1}^n K_i+\frac{\sigma^2}{\mu}$ in both ensembles. If the constraint does not prevent $\Exp_{\mathrm{mic}}[n^{-1}\sum_{i=1}^n K_i + \tfrac{\sigma^2}{\mu}]$ to take the value $(n-1)p+(1-p)$, then we have the same result as in Theorem~\ref{BEEspec}(2), which supports the working hypothesis put forward in Section~\ref{sec:intro}. Indeed, as shown in Section \ref{sec:method}, $S_n(\Pmic\parallel\Pcan)=o(n)$ is the condition for EE. Instead, in the sparse regime we have to rely on events that hold with $(\xi,\nu)$-high probability, where $\xi$ is in principle allowed to vary with $n$, and the condition has to be checked for the specific value of $p(n)$}. 
\end{remark}




\end{document}